\documentclass[12pt]{amsart}
\usepackage{amssymb}
\usepackage{fullpage}
\usepackage{url}
\usepackage[pdftex]{color,graphicx}

\def\F{\mathbb{F}}

\newcommand\Aut{\mbox{\rm Aut}}

\newcommand\CP{\mbox{\rm CP}}

\newcommand\sz[1]{\left|#1\right|}

\newcommand\Perf[1]{\mbox{Perf}(\le #1)}

\def\totalcnt{$15768$}

\author{Alexander Hulpke}
\address{Department of Mathematics, Colorado State University, 1874 Campus
Delivery, Fort Collins, CO, 80523-1874}
\email{hulpke@colostate.edu}
\title{The perfect groups of order up to two million}

\begin{document}
\maketitle

\begin{abstract}
We enumerate the {\totalcnt} perfect groups of order up to $2\cdot 10^6$, up to
isomorphism, thus also completing the missing cases in~\cite{holtplesken89}.
The work supplements the by now well-understood computer classifications of solvable
groups, illustrating scope and feasibility of the enumeration process for
nonsolvable groups.
\end{abstract}

The algorithmic setup for constructing finite groups of a given order, up to
isomorphism,
has been well-established, both in theory and in practice, for the
construction of
groups~\cite{bescheeickobrien02,eickhornhulpke}. It proceeds
inductively, by constructing extensions of known groups of smaller orders
and eliminating isomorphic candidates when they arise. Due to limitations in
implementations of underlying routines, this however had been done so far
mostly for solvable groups.

The aim of this paper is to show the feasibility of generalizing this
approach to the case of nonsolvable groups. 
Instrumental in this has been the calculation of 2-cohomology through confluent
rewriting systems, generalizing the method \cite[\S 8.7.2]{holtbook} for
solvable groups that uses a PC presentation.

The construction process
is illustrated by revisiting the enumeration of
perfect groups that was started in~\cite{holtplesken89} and to extend it to
order $2\cdot 10^6$. In total we find {\totalcnt} perfect groups, seeded from
the 66 nonabelian simple groups of order up to $2\cdot 10^6$.

Compared with~\cite{holtplesken89}, this newly provides explicit lists of the
groups of orders 61440, 86016, 122880, 172032, 245760, 344064, 368640,
491520, 688128, 737280, 983040 that were omitted in their classification of
groups of 
order up to $10^6$. In this range, the calculations also found five groups
(in addition to two groups 
found already in 2005 by Jack Schmidt) 
that
had been overlooked in~\cite{holtplesken89}. Besides serving as examples for
testing conjectures, such lists of groups are used as seed in
algorithms for the calculation of subgroups of a given finite
group~\cite{neubueser60,coxcannonholtlatt,hulpketfsub}, or indeed for the
construction of all groups of a given order.

All calculations were done using the system {\sf GAP}~\cite{GAP4}, which also
serves as repository of the resulting group data.
The program that performed the
classification is available at~\url{https://github.com/hulpke/perfect} and
should allow for easy generalization or extension.

In addition to the actual classification result, this work also serves as a prototype of
enumeration of nonsolvable groups, extending the work
of~\cite{bescheeickobrien02} to the nonsolvable case. It illustrates the
feasibility range of current implementations of underlying routines for
cohomology, extensions, and isomorphism tests, with a number of
general-purpose improvements in the system GAP~\cite{GAP4} (that will be part
of the 4.12 release) by the author having been motivated by this work.

Indeed, the fact, that it took over 30 years since the publication
of~\cite{holtplesken89} to complete the classification of perfect groups up
to order one million, indicates the broad infrastructural requirements of such
classifications, with isomorphism tests~\cite{holtcannonautgroup} being the
most prominent utility (and ultimately the bottleneck of any classification).

\section{The construction process}

We first briefly summarize the construction
process for perfect groups of a given order $n>1$.
This process closely follows the description
in~\cite[\S 11.3]{holtbook} (and, apart from seeding with nonabelian
simple groups, is fundamentally the same strategy as used 
in~\cite{bescheeickobrien02} for solvable groups).

The construction of perfect groups for a chosen order $n$ consists of two
parts, depending on whether the resulting groups have a solvable normal
subgroup or not.

\subsection{Fitting-free groups}

Groups that have no solvable normal subgroup are called {\em Fitting-free}.
Such a group $G$ embeds into the automorphism group of its socle
$S\lhd G$, which in turn is a direct product of simple nonabelian groups. The
conjugation action of $G$ on the $k$ direct factors of $S$ induces a
permutation representation of $G$ of degree $k$. For its image to be
nontrivial perfect, we would need $k\ge 5$ (and thus $n=\sz{G}\ge 60^6$).
For the order range considered, this means that this image is trivial,
thus all direct factors of $S$ must be normal in $G$.
But then $G/S$ is isomorphic to a subgroup of the direct
product of the automorphism groups of the simple nonabelian socle factors.
Such a factor group is solvable (by the Schreier conjecture), showing that
we must have $G=S$ as a direct product of simple nonabelian groups.

For $n\le 2\cdot 10^6$, 
the possible direct factors to consider are:
\begin{eqnarray*}
&&A_5, A_6, A_7, \mbox{PSL}(3,3), \mbox{PSU}(3,3), M_{11}, A_8, \mbox{PSL}(3,4),
\mbox{PSp}(4,3), \mbox{Sz}(8),\\
&&\mbox{PSU}(3,4), M_{12},
\mbox{PSU}(3,5), J_1, A_9, \mbox{PSL}(3,5), M_{22}, J_2, \mbox{PSp}(4,4),
A_{10}, \mbox{PSL}_3(7).
\end{eqnarray*}
and $\mbox{PSL}(2,q)$ for prime powers $7\le q\le 157$, $q\not=128$.  (Note that
$\mbox{PSL}(2,4)\cong\mbox{PSL(2,5)}$ and $\mbox{PSL}(2,9)\cong A_6$.)

\subsection{Inductive construction}

Groups of order $n$ that possess a solvable normal subgroup can be
constructed as extension of groups of smaller order $d\mid n$ by a simple
module of order $p^a=n/d$. As factor groups of perfect groups these smaller
groups need to be perfect themselves. We thus assume that, by induction,
all perfect groups of order dividing $n$ are known.
(Of course the existence of perfect groups of order $d$ is only necessary,
but not sufficient, for the existence of perfect groups of order $n=p^a\cdot
d$.)

We also can assume that $p\mid d$ if $a=1$, since 
the action of a perfect group on a 1-dimensional module must be trivial, and
any extension for $p=n/d$ and $p$ coprime to $d$ thus would be a
direct product and thus not perfect.

This gives the following construction process:

\begin{enumerate}
\item
Iterate over all proper divisors $d\mid n$ 
with $n/d=p^a$, such that $a>1$ or $p\mid d$.
Then iterate over all perfect groups $F$ of order $d$:
\item 
Classify the irreducible $a$-dimensional $F$-modules $M$ over $\F_p$.  For
this, we use the Burnside-Brauer theorem, as described
in~\cite[\S 7.5.5]{holtbook}, to classify all modules, and eliminate those of
the wrong dimension. (Clearly it is sufficient to consider modules for the
factor group $F/O_p(F)$ by the largest normal $p$-subgroup. The index of
the kernel of the module action is further bounded by $\sz{\mbox{GL}_a(p)}$,
which can eliminate some small dimensions $a>1$ for groups that have no small
proper factors.)
\item For each such module $M$ calculate the 2-cohomology group $H^2(F,M)$
(Section~\ref{cohomology}). 
\item For each cocycle $\zeta$, representing an element of $H^2(F,M)$,
construct the corresponding extension $E$.
\item
Test each $E$ for being perfect (and discard $E$, if not).
\item
Eliminate isomorphic groups.
\end{enumerate}

Compared with~\cite{holtplesken89}, this construction does not aim to
construct groups as subdirect products, if possible. Nor
does it use further theoretical results to restrict potential orders,
primes, or module types.  Instead, the generic construction algorithm is
used throughout, to have as much of an independent verification of the prior
results as possible,
and to reduce the potential of missing cases in a lengthy case distinction.
(See Section~\ref{compare} for evidence of this risk.) 

The main difficulty, as with any enumerative construction, however remains
the elimination of isomorphic groups. (This is stated already
in~\cite{holtplesken89} as having been the main obstacle.) Since isomorphism
tests are expensive, we utilize a number of techniques to reduce the number
of tests needed.  These techniques include the incorporation of isomorphism
elimination already in the construction process (Compatible pairs,
Section~\ref{comppairs}), testing for isomorphism invariants that can be
computed cheaply and can decide non-isomorphism in many cases (Fingerprints,
Section~\ref{fingerprints}), and in particular by selecting a
``canonical'' construction
path for groups  (Section~\ref{cancon}). 

Since many of these concepts have been described before in other contexts, the
following descriptions will be brief and focus on the case of perfect groups
and the actual implementation used.

\subsection{Calculating Cohomology}
\label{cohomology}

Assume that a factor group $F$, and an irreducible $\F_pG$-module $M$ have
been chosen. The first step of the calculation is the determination of
$H^2(F,M)$, using the confluence condition for rewriting systems.
Such an approach has been suggested, e.g.,  in~\cite{schmidt10}. A
detailed description can be found in \cite[\S 7]{dietrichhulpke21}.

First, we determine a confluent rewriting system for $F$. Such a rewriting
system can be composed (with a wreath product ordering) from rewriting
systems for the composition factors of
$F$, which in turn can be found using the methods of \cite{schmidt10}. We also found
it useful, in the case of large primes, to introduce powers of generators
as extra generators, e.g. writing $x^{50}$ as $y^4x^2$ with
$y=x^{12}$, as doing so keeps word lengths shorter.

The rewriting system for $F$ can be extended to one for an extension of
$M$ by $F$ by adding generators for $M$, adding relations describing the
$F$-module structure of $M$ (that is $m\cdot f\to f\cdot m^f$), and by modifying each rule $l\to r$ for $F$ by
a variable tail $t\in M$ to $l\to r\cdot t$.

The ``critical pair'' confluence conditions of the Knuth-Bendix algorithm (which hold
for the rewriting system in $F$, since it was assumed to be confluent)
then translate into linear conditions on the the tails.
If we represent 2-cochains as vectors, composed from sequences of tail
values, these conditions form a homogeneous system of linear
equations, whose solution space represents the $2$-cocycles $Z^2(F,M)$.

We represent the the $2$-coboundaries $B^2(F,M)$ by a subspace therein.
A generating
set for this subspace is obtained by all combinations of  changing coset
representatives for a generator $x$ of $F$ by basis vectors of $M$.

Then $H^2(F,M)$ can be represented by a complement subspace to $B^2$ in
$Z^2$, each of its elements
representing a choice of tails to yield a rewriting system (and thus a
presentation) for the associated extension of $M$ by $F$.

\subsection{Compatible pairs}
\label{comppairs}

Isomorphisms that respect the extension structure induce an action on 
the cohomology group through the group of compatible pairs:

An  element $(\kappa,\nu)\in\Aut(F)\times\Aut(M)$ (here $\Aut(M)$ is the
group of vector space automorphisms, i.e. $\mbox{GL}(M)$)
is called a {\em compatible
pair}~\cite{robinson81}, if its action is compatible with the module action of $F$ on $M$,
that is, if for any $f\in F$ and for any $m\in M$ we have that 
\begin{equation}
\nu(m^f)=\nu(m)^{\kappa(f)}.
\label{compp}
\end{equation}
The compatible pairs form a group $\CP$. It acts on $H^2(F,M)$, with extension
corresponding to cocycles in the same orbit being
isomorphic~\cite[\S 8.9]{holtbook}. It thus is
sufficient to consider only orbit representatives under this action for the
construction of groups.

We compute the group of compatible pairs as follows:
If $\kappa=1$, the condition~(\ref{compp}) is that
$\nu\in\Aut_{\F_p F}(M)$ must be an $\F_pF$-module automorphism. This group 
$\Aut_{\F_p F}(M)$ of module automorphism can be computed
using MeatAxe-style techniques~\cite{michaelsmithphd}.

Furthermore, if two compatible pairs $(\kappa,\mu)$ and $(\kappa,\nu)$ agree
in their first component, their quotient $(1,\mu/\nu)$ will be a compatible
pair. This means that if will be sufficient to identify the nontrivial
$\kappa$ that can be first component of a compatible pair $(\kappa,\nu)$,
and if so to determine a suitable $\nu$. For this, 
given $\kappa\in\Aut(F)$, we consider the vector space $M$ in a
second way as an $\F_pF$-module (which we shall call $M_\kappa$), namely
with $F$ acting through its image under $\kappa$. Then $(\kappa,\nu)$ is a
compatible pair, if and only if $\nu$ is an $\F_pF$-module isomorphism
between these $F$ modules $M$ and $M_\kappa$. Again, MeatAxe-style
techniques~\cite{michaelsmithphd} can test for such module isomorphisms constructively, producing a suitable $\nu$.

We then use this test in a backtrack search through $\Aut(F)$, using a
faithful permutation representation of $\Aut(F)$
\cite[\S4.6]{holtbook}. This search determines the image of the projection
of the compatible pairs on the first component. For each $\kappa\in\Aut(F)$
tested successfully, we store the compatible pair $(\kappa,\nu)$ for $\nu$
produced by the isomorphism test.

These compatible pairs, together with pairs $(1,\nu)$ for $\nu$ from a
generating set for the group of module automorphisms, generate the group of
compatible pairs.
\smallskip

The determination of orbit representatives of $\CP$ on $H^2(F,M)$ then is
straightforward, as in the cases encountered the dimension of the cohomology
group is sufficiently small to list its elements.

\subsection{Permutation Representations}

The cocycles representing $H^2$ each are collections of values of tails for
rewriting rules for the corresponding extension $E$. We thus get each
extension $E$ in the form of a finitely presented group.

We first test each group $E$ obtained this way for it being perfect, and
discard those groups that are not.
This test can be performed
effectively, given only a presentation of $E$, using a Smith Normal form
computation. 
\smallskip

Any further checks however will be done far
more efficiently, if the group can be described by a faithful permutation
representation of moderate degree.
(Such a representation also is desirable for a resulting library of
perfect groups.)
To find such a representation we follow the strategy in~\cite{hulpkequotgp}:

In a slight abuse of notation we write $M\lhd E$ to denote the normal
subgroup associated to the construction. We are searching through subgroups
of $F$ of small index~\cite{cannonholtslatterysteel}, in search of a
subgroup $U\le F$, such that its pre-image $T\le E$ has an abelianization
$T/T'$ that is larger than $U/U'$. (The abelianization of subgroups of
finitely presented groups, indeed an epimorphism $T\to T/T'$,
can be found by standard methods~\cite{simsbook}.) 
In this case, there will be a permutation representation $\lambda$ of $T$
with $T'\le\ker\lambda$ (that is, we can find $\lambda$ from the
representation $T\to T/T'$) with $M\not\le\ker\lambda$. The subdirect product
of the induced representation $\lambda\uparrow^E$, with a representation
$E\to F$, then is a representation of $E$ with a kernel strictly smaller
than $M$. As $M$ is a simple module, the subdirect product of
$\lambda\uparrow^E$ with a faithful representation of $F$ will be a
faithful representation of $E$.

We found that this process resulted in reasonable degree permutation
representations for the purpose of isomorphism tests. When storing the
groups with permutation representations in a library, we however also ran
degree reduction heuristics (the {\sf GAP}
command \verb+SmallerDegreePermutationRepresentation+) on the resulting
groups to try to get a smaller degree permutation representation.

For about 0.5\% of the groups this resulted in permutation representations
of unduly larger degree, the worst being degree $36864$ for four groups of order $1843200$.
To make the use of the groups less costly, we re-ran the degree reduction for a
longer time for these outlier groups,
also testing the action on cosets of a number of randomly
generated subgroups. In five remaining stubborn cases, finally, partial
subgroup lattice information was used to determine a smaller degree permutation
representation by hand. These reductions resulted in a maximal degree of $8448$
(occurring for one group of order $1966080$).
We did not aim to establish minimal faithful permutation degrees for all
groups, as this would have come at significant higher cost --
essentially the computation of subgroup lattices.

\begin{figure}
\begin{center}
\includegraphics[width=70mm]{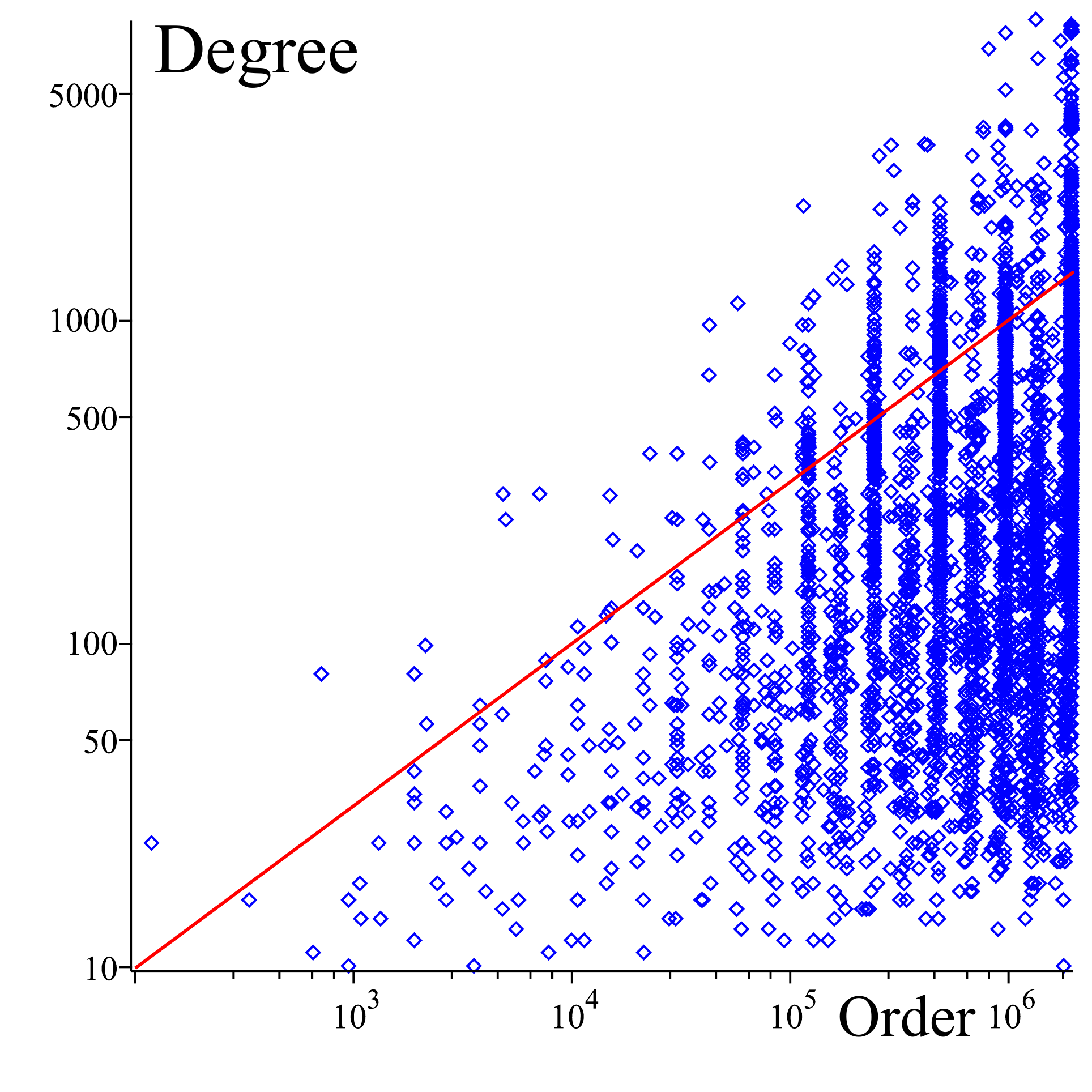}
\qquad
\includegraphics[width=70mm]{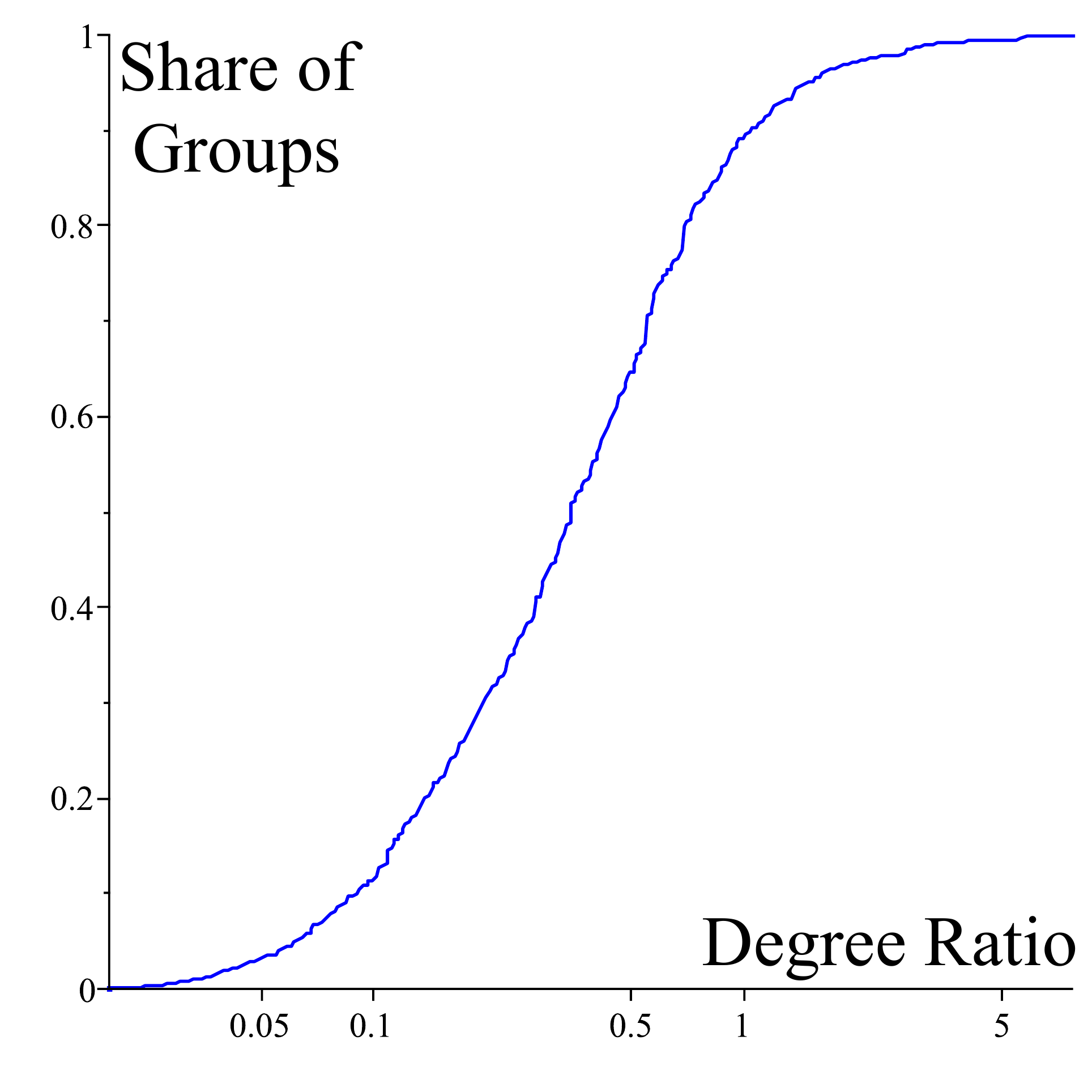}
\end{center}
\caption{Permutation degrees in relation to $\sqrt{\sz{G}}$.}
\label{figdeg}
\end{figure}

Figure~\ref{figdeg} shows the overall statistics of the resulting
permutation degrees.  The log-log plot on the left hand side displays a blue
scatter plot of degrees versus group orders. Its shape, and the
least-squares fit $\sz{G}^{0.42}$, suggest to consider the degrees in
relation to $\sqrt{\sz{G}}$, which is given by the red line.

The right hand side of Figure~\ref{figdeg} shows, for a ratio $r$ given on
the $x$-axis, the share of groups for which
$\mbox{degree}(G)/\sqrt{\sz{G}}\le r$.  It indicates that about 80\% of the
groups have a permutation degree between $0.1\sqrt{\sz{G}}$ and
$\sqrt{\sz{G}}$, while all have a degree bounded by $10\sqrt{\sz{G}}$.
This justifies that for the groups constructed a degree of magnitude
$\sqrt{\sz{G}}$ can be expected.

\section{Isomorphism Rejection}
\label{isomorej}

At the heart of any combinatorial classification algorithm is the question
of isomorphism rejection. This is needed here, since 
the same abstract group might be constructed in multiple ways, e.g.
through different normal subgroups, and such cases will not be identified by
compatible pairs.

Isomorphism tests however are expensive,
which makes it imperative to reduce the overall number of tests required.
We do so through methods that are common to combinatorial enumerations

\subsection{Canonical construction}
\label{cancon}

We associate to each group $E$, constructed as extension of $M\lhd E$ by
$E/M\cong F$, a ``canonical'' construction path.  This allows us to discard
any group we have constructed, if it turns out that it had been constructed
in a non-canonical way.

We define this canonical construction in the following way:
\begin{enumerate}
\item 
Amongst all minimal normal subgroups $N\lhd E$, $M$ must be of minimal order.
\item
Amongst all minimal normal subgroups $N\lhd E$ with $\sz{N}=\sz{M}$, the
isomorphism type of $F=E/M$ is minimal (in a list of isomorphism types of
perfect groups of order $d=\sz{F}$).
\end{enumerate}
These conditions imply that a group $E$, constructed from $F$ and $M$, gets
discarded if there is $N\lhd E$ with $\sz{N}\lneqq \sz{M}$ or with
$\sz{N}=\sz{M}$ and a factor group $F/N$ of smaller (in a given list)
isomorphism type. Such an isomorphism type does not need to be determined
exactly, it is sufficient to know that any possible isomorphism type is
smaller (or larger) than the type of $F$.

Once it has been established that a group has been constructed canonically,
the only isomorphisms that we need to test for must be amongst groups constructed
for the same $F$ and for modules of the same order. Furthermore, such an
isomorphism may not fix $M$, since it otherwise would result in a compatible
pair. This means that the group $E$ must have a second minimal normal
subgroup $M_2$, such that $E/M\cong E/M_2$, and that the $E/M$-modules $M$
is isomorphic to the $E/M_2$-module $M_2$.
Explicit isomorphism tests thus only need to be run
amongst the extensions that satisfy this
property, which means in particular that they must have been 
constructed from the same factor group $F\cong E/M$.

This restriction significantly reduces the number of isomorphism tests required.

\subsection{Fingerprints}
\label{fingerprints}

The second of the ``canonicity'' conditions requires us to identify the isomorphism type of
factor groups. Since we assume that a list of all perfect groups of order
$d\mid n$ is known a priori, we can usually identify this type, without (or with
minimal) need for isomorphism tests, by calculating a ``fingerprint'' of the
group, composed of values obtained from isomorphism invariants. Groups with
different fingerprints cannot be isomorphic, and the hope is that
fingerprint information identifies a group uniquely in a given list.
(Typically this hope is not
fully satisfied, and would be provable only post-factum for an explicitly
constructed list of isomorphism types.)
In fact, even partial fingerprinting
information can be sufficient, as long as it proves that all possible isomorphism
types for a given factor group of $E$  are smaller (respectively larger)
than the type of $F$ used in the construction.

The same kind of fingerprint information can be used to disprove
isomorphism, before attempting an actual search for
isomorphisms~\cite{holtcannonautgroup} amongst groups constructed for the
same $F$ and for modules of the same order.
\smallskip

In the construction of perfect groups of order up to $2\cdot 10^6$, the
following fingerprint properties were used. This list was determined
experimentally as providing a decent tradeoff between avoiding isomorphism
tests and not being too costly on their own. They however do not guarantee
unique identification, and in a handful of cases explicit isomorphism tests
have been required to identify factor group types uniquely.
(If (sub)groups are given as
identification, the actual 
fingerprint information for each subgroup consists of its order, whether it
is perfect, as well as its
identification in the small groups library~\cite{bescheeickobrien02} if the
order is sufficiently small.)
Clearly one can stop computing fingerprint information, as soon a group is
uniquely identified in a given list.
\begin{itemize}
\item Conjugacy class representatives and centralizers.
\item Normal subgroups and their centralizers.
\item Maximal, and low index~\cite{cannonholtslatterysteel} subgroups.
\item Automorphism groups.
\item Characteristic subgroups.
\item Character table (i.e. classes/characters being permutable to make the
resulting tables identical. This is done by a backtrack search for a
class permutation whose
effect on the character table can be undone by a character permutation and
is provided by the {\sf GAP} command
\verb+TransformingPermutationsCharacterTables+).
We attempt this test (that requires the calculation of character tables in
the first place)
for factor groups that have at most $200$ classes,
respectively full groups (before an isomorphism test) that have at most
$500$ classes, and abort such a backtrack search for a suitable class
permutation if its runtime exceeds 5-10 minutes, as at that point an explicit
isomorphism test becomes competitive for the groups considered here.
(Again, these parameters were found experimentally as
providing a reasonable trade-off between computation cost and potential
speedup.)
\end{itemize}

We note as an aside that character tables will not give a perfect identification for the
pool of groups considered. The smallest example for this happening are two
groups (Numbers 37 and 38 in the published list) of order 61440 and
structure $2^2.2^4.2^4.A_5$, which form a Brauer pair. (These two groups
differ in the number of classes of subgroups -- 15731 versus 15715 -- but this
invariant count turns out to be more expensive to compute than an isomorphism
test, which is the reason that we do not use it in the classification.)

\section{Comparisons with~\cite{holtplesken89}}
\label{compare}

The implementation of the algorithm was tested by revisiting the
constructions for the orders considered in~\cite{holtplesken89}.

In checking results, we found five new groups, of order 243000 (two groups),
729000, 871200, and 878460 that had been overlooked in the earlier
classification. (All of them are ``obvious'' constructions without any need
of isomorphism testing, and their omission is unlikely to be the result of a
calculation error, but seems to stem from sub-cases of a complicated
construction having been overlooked accidentally. That is, they do not
indicate conceptual errors in the earlier classification, but are simply
clerical errors.) We also confirmed the two missing groups (or order
$450000$ and $962280$) found in 2005 by Jack Schmidt.

All other results of~\cite{holtplesken89} were confirmed. Given the
substantial progress in both hardware and software since then, this makes
the results of~\cite{holtplesken89} even more impressive and should
be considered as a resounding validation of that work.
The discovery of a small number of overlooked groups however also serves as
corroboration of the more basic construction approach used here, 
in which we do not aim to
use special methods for particular cases -- such reductions do not apply to the
hardest cases, and thus have little impact on the overall time required, but
they increase the risk of accidental omissions.
\smallskip

As for the orders $<10^6$, for which~\cite{holtplesken89} omitted lists,
we found counts of groups as given in
table~\ref{bigres} (the underlined numbers
already were determined in~\cite{holtplesken89}, though not the actual groups).
The calculations for these 11 orders took about 2
(single-core) weeks on a Dell R740, Xeon
Gold 6132, 2.60 GHz Processor (Geekbench 5 Single-Core 910) and used about
20GB of memory.

\begin{table}
\begin{tabular}{l|r|r|r|r|r|r}
Order&61440&86016&122880&172042&245760&344064\\
\hline
Count&98&\underline{52}&258&154&582&291\\
\hline
\hline
Order&368640&491520&688128&737280&983040\\
\hline
Count&\underline{46}&1004&508&\underline{54}&1880\\
\end{tabular}
\caption{Newly found perfect groups of order $\le 10^6$, counts}
\label{bigres}
\end{table}

Lists of the groups 
will be made available as part of the
system {\sf GAP}~\cite{GAP4}, release 4.12.

\section{Results}

We furthermore enumerated the perfect groups of orders between $10^6$ and
$2\cdot 10^6$.  Their counts are given
\footnote{Due to an spreadsheet editing error, an earlier version
of this note gave some wrong values in the table and subsequently a wrong
total.}
in table~\ref{perfcnt}, giving in
total {\totalcnt} perfect groups of order up to $2\cdot 10^6$.

The total calculation time for these groups 
has been about 11 (single-core,
on the same machine as above) weeks, with over 90\% of the time taken by by
order $1966080$.
Most of this time was (unsurprisingly) taken by work for isomorphism
rejection, namely in
identifying the type of factor groups in the test for
the canonical construction described in Section~\ref{cancon}, while only
comparatively few explicit isomorphism tests have been necessary.

The largest cohomology groups encountered were of dimension $12$, with
up to $256$ orbits remaining under the action on compatible pairs.
(This occurs for example for the trivial module of group $957$ of order
$983040$). The increase in the number of orbits on the cohomology groups,
and the time required indicates that any significant extension of the lists
would require a magnitude more of resources, which ultimately is the reason
for stopping the classification at order $2\cdot 10^6$. 

With most of the groups being distributed amongst a small number of orders
-- primarily products of $60$ of $168$ with powers of $2$, and with most groups
being extensions of groups $F$ of order $n/2$ by trivial modules $M$ of 
order $2$ -- it of course
would be easy for an interested researcher to calculate the groups of other,
larger, orders with the provided program, as long as their total number was
small.
\medskip

\begin{figure}
\includegraphics[width=75mm]{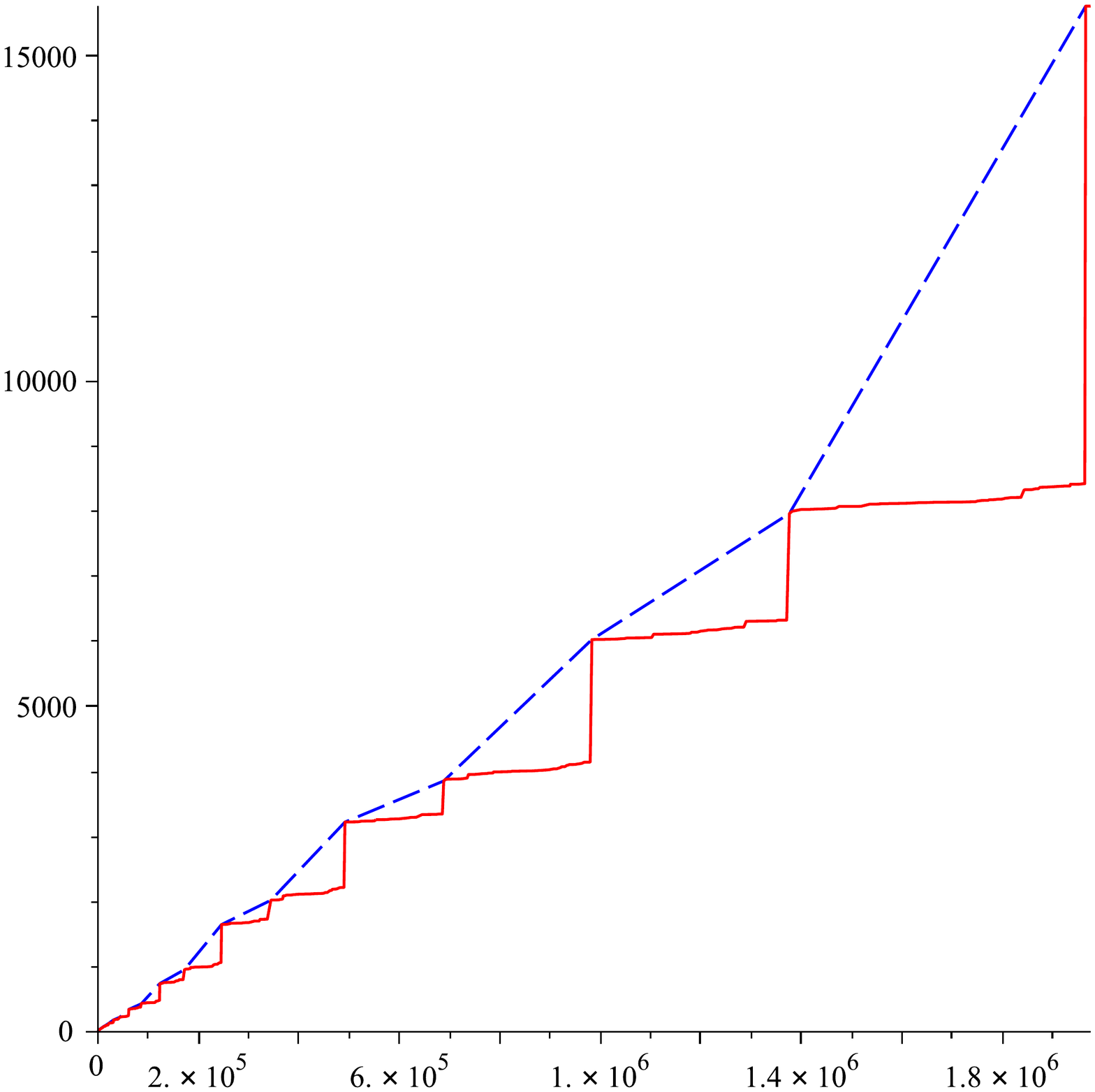}
\includegraphics[width=75mm]{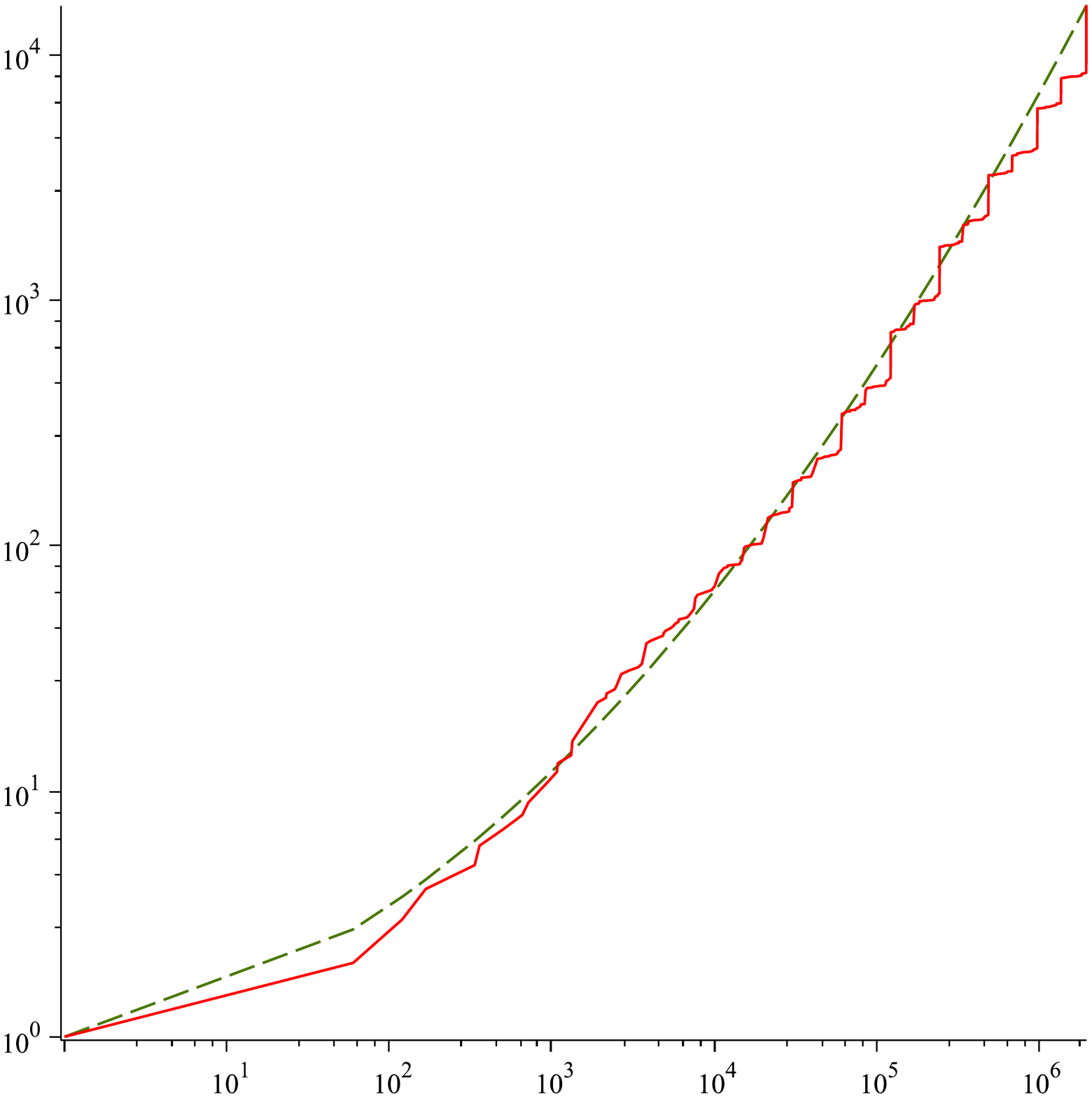}
\caption{Counts of perfect groups, $\Perf{n}$}
\label{figplot}
\end{figure}

We also tried to compare the experimental counts of groups with theoretical
predictions in~\cite{holtperfcount}.
For this, denote by $\Perf{n}$ the number
of perfect groups of order $\le n$. The estimates in
\cite[Theorem A]{holtperfcount},
\[
n^{\log_2(n)^2/108-c\cdot\log_2(n)}\le \Perf{n}\le
n^{\log_2(n)^2/48+\log_2(n)}
\]
give for $n=2\cdot 10^6$ only that 
\[
1.8\cdot 10^{-15}\le \Perf{2\cdot 10^6}\le 2.6\cdot 10^{189},
\]
and thus agree with the result, but do not provide a meaningful validation.
The exceedingly tiny lower bound here stems from choosing the parameter
$c=11/36$, which is indicated in~\cite{holtperfcount} as a
correct, but unlikely to be optimal, estimate.

Figure~\ref{figplot}, left, shows a plot of
$\Perf{n}$ (red), as well as an ``envelope'' (blue, dashed) that accounts
for the main growth contributors from a few orders. The changing growth
of this envelope is due to different series of extensions with factor groups
$A_5$, respectively $L_3(2)$, and the fact that the simple groups do not
have irreducible modules in every dimension. (E.g. 
perfect groups of order $n$ with radical factor $A_5$ give rise to
extensions of order $2n$ and $16n$, but not $4n$ or $8n$.)

Theorem $A$ in~\cite{holtperfcount} indicates that $\Perf{n}$
behaves asymptotically similar to a function $n^{p(\log_2(n))}$ with $p$ a
polynomial of degree $2$, thus $\log_2(\Perf{n})$ should be a polynomial
of degree $3$ in $\log_2(n)$. A least squares fit of the values for the
envelope (to account for the growth mainly occurring in a few orders -- we
also attempted a fit for all data and got similar results with a worse
match)
of $\Perf{n}$ finds a degree $3$ best fit polynomial
\[
-0.00060 x^3 +0.04547 x^2 -0.03027 x -0.07598
\]
whose coefficients (in magnitude, as well as the negative coefficient of the
leading term) seem to indicate an ill fit. Tries of several polynomials
make a quadratic fit
$0.02807 x^2+0.08021x$ seem the best match.
Figure~\ref{figplot}, right, gives a log-log plot of $\Perf{n}$ (red) and
the (green, dashed) fit curve
\[
n^{0.02807\log_2(n)+0.08021}.
\]
This discrepancy to the asymptotic results~\cite{holtperfcount} in the degree in the exponent is
likely due to the orders considered being still too small to be
representative of the asymptotic behavior, and the logarithmic curve fitting
giving more weight to the group of small order.

\begin{table}
\begin{tabular}{c|c||c|c||c|c||c|c}%
Order&Count&Order&Count&Order&Count&Order&Count\\
\hline
1008000&1&1233792&1&1467648&2&1774080&9\\
1008420&1&1244160&15&1468800&1&1785600&3\\
1016064&1&1253376&4&1474560&26&1787460&1\\
1020096&1&1260000&2&1512000&1&1788864&4\\
1024128&1&1267200&15&1518480&1&1800000&3\\
1030200&1&1270080&2&1536000&33&1806336&13\\
1036800&3&1277760&2&1548288&1&1814400&6\\
1044480&4&1285608&1&1555200&3&1815000&3\\
1048320&2&1290240&88&1572480&1&1822500&2\\
1053696&9&1294920&1&1574640&4&1837080&1\\
1080000&1&1296000&1&1592136&1&1843200&113\\
1083000&1&1310400&1&1614720&3&1843968&3\\
1088640&1&1330560&2&1615680&1&1845120&3\\
1092624&1&1342740&1&1632960&7&1858560&1\\
1100736&1&1350000&1&1645056&1&1866240&13\\
1102248&3&1351680&8&1651104&1&1872000&1\\
1105920&49&1354752&3&1653900&1&1875000&22\\
1123980&1&1370880&1&1658880&2&1876896&1\\
1125000&1&1376256&1639&1663200&1&1886976&1\\
1149120&3&1382400&38&1693440&2&1920000&15\\
1166400&4&1386240&3&1713660&1&1924560&2\\
1176120&3&1399680&21&1721400&1&1934868&1\\
1179360&4&1414944&1&1723680&1&1935360&26\\
1180980&14&1425600&3&1728000&1&1953000&1\\
1192464&1&1425720&1&1728720&1&1959552&6\\
1200000&17&1441440&3&1742400&2&1964160&1\\
1209600&8&1442784&1&1747200&1&1966080&7344\\
1215000&9&1451520&3&1749600&8&1975680&1\\
1224120&1&1457280&1&1756920&8&1980000&1\\
1224936&1&1461600&3&1762560&3&&\\
1231200&1&1463340&1&1771440&1&&\\
\end{tabular}
\caption{Counts of perfect groups beyond $10^6$}
\label{perfcnt}
\end{table}

\section*{Acknowledgment}
The author's work has been supported in part by
NSF~Grant~DMS-1720146, which is gratefully acknowledged.

All calculations were performed on a departmental computing server. The
author would like to express his thanks to the Department of Mathematics and
the College of Natural Sciences at Colorado State University for access to
this machine.
He also would like to thank an anonymous referee for careful reading and
helpful comments.

\bibliographystyle{alpha}
\bibliography{mrabbrev,litprom}

\end{document}